\input amstex

\documentstyle{amsppt}

\baselineskip=24pt
\pagewidth{130mm}
\pageheight{175mm}
%\NoBlackBoxes
\NoRunningHeads
\loadeusm

\magnification = 1200

\define\BP{\bigpagebreak}
\define\bs{\backslash}

\define\F{\flushpar}
\define\q{\qquad}
\define\qq{\qquad\qquad}

\define\SP{\smallpagebreak}

\TagsOnRight

\define\A{\;$}

\define\CH{\Cal H}

\define\Con{C_0^\infty}

\define\f{$\,}

\define\gr{\nabla}

\define\Mtt1{{\Tilde {\Tilde M}}_1}
\define\Oa{\Omega}
\define\pa{\partial}
\define\pd#1#2{\dfrac{\partial#1}{\partial#2}}

\define\re{\roman{Re}}
\define\R3{\bold R^3}

\define\RN{\bold R^N}

\define\Sch{Schr\"odinger}

\define\T{\tag}

\define\Eia{1} % 1960 paper
\define\Eib{2} % CMP
\define\Ik{3}  % Eigenfunction expansion
\define\IS{4}  % Limiting absorption method
\define\IU{5} % Ikebe-Uchiyama
\define\Ja{6} % uniqueness
\define\JSa{7} % F Math
\define\JSb{8} % Osaka Math
\define\Ka{9} % Growth property
%\define\KI{ } % Kato-Ikebe
\define\Mi{10} % M\"uller
\define\Rl{11} % Rellich
\define\RZ{12} % Roach-Zhang
\define\YS{13} % Springer 727
%\define\Sa{ } % Pacific J.
\define\We{14} % Weidmann
\define\Wi{15} % Wienholtz

  % This is a version of "unique1" to comply with the referee's suggestion.
  % This was originally finished on February, 1997 at Heidelberg.
  % This revision was done in October, 1997.

\topmatter
\title  The uniqueness of the solution of the Schr\"odinger equation
              with discontinuous coefficients         \endtitle
\author Willi J\"ager$^1$ and Yoshimi Sait\B o$^2$           \endauthor
\affil  $^1$ University of Heidelberg \\
        Institute for Applied Mathematics \\
        D-69120 Heidelberg  \\
        Germany \\ \ \ \ \\
        and \\
        \ \ \ \\
        $^2$ Department of Mathematics \\
        University of Alabama at Birmingham \\
        Birmingham, Alabama 35294 \\
        U. S. A. \\                  \endaffil

\endtopmatter

\newpage

\document

\head     {\bf \S1. Introduction}      \endhead
\BP
        Let us consider the reduced wave equation
$$
           -\,\Delta u(x) + q(x)u(x) = 0                        \T 1.1
$$
     on the domain \f \Omega \A such that
$$
     \Omega \supset E_{R_0} = \{ x \in \RN \ : \ |x| > R_0 \},   \T 1.2
$$
     where \f R_0 > 0 \A and \f N \ge 2 $. Suppose that \f q(x) \A has the
     form
$$
          q(x) =  - \ell(x) + s(x),                            \T 1.3
$$
     where \f \ell(x) \A a positive function, and \f |s(x)| \A is supposed
     to be dominated by \f \ell(x) $. The equation (1.1) has been studied
     extensively especially in the relation to the operator
$$
        H_1 = - \Delta + V(x)                                   \T 1.4
$$
     in \f L_2(\Oa) $, or
$$
        H_2 = - \frac1{\mu(x)}\Delta                           \T 1.5
$$
     in the weighted Hilbert space \f L_2(\Oa; \, \mu(x)dx) \A with
     boundary conditions on the boundary \f \pa\Oa \A and at infinity.
     In this work we are concerned with the asymptotic behavior of the
     solution \f u \A of the equation (1.1) at infinity. One of the important
     conclusions of the study is that we can establish the nonexistence
     of a class of (nontrivial) solutions of (1.1) which includes the
     \f L_2 $-solutions. And this result plays an important role in the
     attempt (the limiting absorption method, see, e.g., [\Eia], [\IS]) to
     prove the existence of the boundary value of the resolvent
     \f (H_1 - z)^{-1} \A or \f (H_2 - z)^{-1} \A of the operator
     \f H_1 \A or \f H_2 \A when the complex parameter \f z \A approaches
     the real axis.
     
        Consider the equation
$$
          -\,\Delta u(x) + (- k^2 + s(x))u(x) = 0,                \T 1.6
$$
     with \f k > 0 $, i.e., \f \ell(x) = k^2 \A in (1.3). In
     the celebrated work Kato\,[\Ka] he showed, among others, that, under
     the condition
$$
       \tau \equiv (2k)^{-1}\varlimsup_{|x|\to\infty} |x||s(x)| < 1,
                                                              \T 1.7
$$
     a nontrivial solution \f u \A of the equation (1.1) satisfies
$$
       \lim_{r\to\infty} r^{2\tau+\epsilon}
         \int_{|x|=r} \bigg( |\gr u(x)|^2 + |u(x)|^2 \bigg) \, dS = \infty
                                                                 \T 1.8
$$
     for any \f \epsilon > 0 $. One of the important features of the work
     [\Ka] is that the coefficient \f s(x) \A does not need to be
     spherically symmetric which makes the scope of application
     much wider than the preceding works (cf., e.g., M\"uller\,[\Mi], 
     Rellich\,[\Rl]). Another important feature of [\Ka] is that the method 
     is based on differential inequalities satisfied by several functionals 
     of the solution \f u $ so that the problem was successfully treated as 
     a local problem at infinity. As a result we do not need to use any
     boundary conditions at the boundary \f \pa\Omega \A of \f \Omega \A or
     at infinity such as radiation condition (cf., e.g., Wienholtz\,[\Wi]).
     As is well-known this result has many applications. In Ikebe [\Ik],
     in which the spectral theory and scattering theory for the \Sch\
     operator \f - \Delta + V(x) \A in \f \R3 \A was developed under the
     conditin that
$$
       V(x) = O(|x|^{-\gamma}) \qq (|x| \to \infty, \gamma > 2),    \T 1.9
$$
     the result of Kato [\Ka] was used to prove the existence of the
     boundary value of the Green function on the positive real axis
     as well as the nonexistence of the positive eigenvalues. After the
     work [\Ka] various extensions and modifications were presented as
     many efforts were made to treat more general operators in a similar
     method. See e.g., Ikebe-Uchiyama \,[\IU] for \Sch\ operators with
     magnetic potentials, J\"ager \,[\Ja] for the second order elliptic
     operators, Weidmann \,[\We] for the many body \Sch\ operators, and
     Ikebe-Sait\B o \,[\IS], Sait\B o[\YS] for \Sch\ operators with
     long-range potentials.
     
        Now let us consider the case that \f \ell(x) \A is a positive
     function which may be discontinuous. One of the motivations to consider
     such \f \ell(x) \A comes from the study of the reduced wave equations
     in layered media. Consider the equation
$$
       - \mu(x)^{-1}\Delta u - \lambda u = 0  \q (x \in \RN)      \T 1.10
$$
     in layered media, where \f \mu(x) \A is a positive function on
     \f \RN $. Suppose that the function \f \mu(x) \A is a simple function
     with surfaces of discontinuity (separating surfaces) which may extend
     to infinity. Roach and Zhang\,[\RZ] proved the nonexistence of the
     solution of the equation (1.10) under a geometric condition
     (``cone-like'' discontinuity on the separating surface, see also 
     [\Eib]). Then J\"ager and Sait\B o\, [\JSa] proved a similar results 
     under another geometric condition (``cylindrical'' discontinuity) on 
     the separating surfaces. In these works the method is not local at 
     infinity, but some global integral identity of the solution \f u \A are 
     used. And the method seems to need some modifications in the case where
     \f \mu \A is a perturbation such as
$$
        \mu(x) = \mu_0(x) + \mu_{\ell}(x) + \mu_s(x),              \T 1.11
$$
     \f \mu_0 \A being a simple function, and \f \mu_{\ell}(x) \A and
     \f \mu_s(x) \A behaving like a long-range and short-range potentials
     at infinity, respectively (cf. [\JSb]).
     
        In this work we are going to obtain an extension of the result (1.8)
     by Kato\,[\Ka] which can be applied the reduced wave equation (1.10)
     with \f \mu(x) \A satisfying (1.11) as well as the equation (1.6) where
     \f s(x) \A is the sum of a short-range potential and a long-range
     potential. Under the several assumptions (Assumptions 2.1, 4.2, 5.5 and
     5.8) on the coefficient \f q(x) \A the following (Theorem 5.10 in \S 5)
     will be proved:
\BP

        {\it Suppose that a solution \f u \A of equation {\rm (1.1)}
     satisfies
$$
        \varliminf_{r\to\infty}
             \int_{|x|=r} \big\{ \big|\frac{\pa u}{\pa r}\big|^2
                   - \re\,(q(x))|u|^2 \big\}\, dS = 0.        \T 1.12
$$
     Then \f u \A has a compact support.}
\BP
     
        Our method is a local method at infinity which is
     similar to the method of Kato\,[\Ka]. As in [\Ka], some type of
     differential inequalities on functionals of the solution \f u \A will
     play important roles. However, we shall first establish the
     differential inequalities not in the ordinary sense but in the sense
     of distributions, and then they will be interpreted in the ordinary
     sense.
     
        In \S2 we define our reduced wave equation and give the
     main assumption (Assumption 2.1) on the coefficients. In \S3 we
     introduce and evaluate the first functional \f M^+(v, r) $. In order
     to complete the evaluation of \f M^+(v, r) $, another functional
     \f N(v, m, r) \A is introduced and evaluated in \S4. \S5 is devoted for
     proving the main theorem (Theorem 5.9). Some examples are discussed in
     \S6. In \S7 we shall discuss how our result can be applied to some
     reduced wave operators which were studied in [\JSb]. A lemma on 
     distributional derivative is given in Appendix.
\SP

        {\bf Acknowledgement.} This work was finished when the second
     author was visiting the University of Heidelberg for February 1997.
     Here he would like to thank Deutsche Forschungs Gemeinschaft for
     its support through SFB 359.  Also the second author is thankful to
     Professor Willi J\"{a}ger for his kind hospitality during this period.

$$
\ \ \ \
$$

\head     {\bf \S2. \Sch-type homogeneous equation}      \endhead
\BP

        Consider the homogeneous \Sch\ equation
$$
           -\,\Delta u(x) + q(x)u(x) = 0 \qq
                                (x \in E_{R_0}),              \T 2.1
$$
        where \f R_0 > 0 $, and
$$
          E_{R} = \{ x \in \RN : \ |x| > R \}.               \T 2.2
$$
     Let \f S^{N-1} \A be the unit sphere of \f \RN $. We set
     \f X = L_2(S^{N-1}) \A and the inner product and norm of \f X \A is
     denoted by \f ( \ , \ ) \A and \f | \ | $, respectively.
\BP

     {\bf Assumption 2.1.} \ (i) Let \f N \A be an integer such that
     \f N \ge 2 $. Let \f u \in H^2(E_{R_0})_{\roman loc} $, \f R_0 > 0 $,
     be a solution of the equation (2.1), where \f q(x) \A is a
     complex-valued, measurable, locally bounded function on \f E_{R_0} $.
\SP

        (ii) Set
$$
           Q(x) = q(x) + \frac{(N-1)(N-3)}{4r^2}.             \T 2.3
$$
\SP

\roster
\item"{(ii-a)}" Then \f Q(x) \A is decomposed as
$$
         Q(x) = Q_0(x) + Q_1(x),                             \T 2.4
$$
     where \f Q_0(x) \A is a real-valued, measurable, locally bounded
     function on \f E_{R_0} \A \linebreak such that 
$$
           Q_0(x) \le 0,                                     \T 2.5
$$      
     and \f Q_1(x) \A is a complex-valued, measurable, locally bounded
     function on \f E_{R_0} $.
\SP

     \item"{(ii-b)}" For any \f x \in X = L_2(S^{N-1}) $, 
     \f (Q_0(r\cdot)x, x) \A has the right limit for all \f r > R_0 \A
     as a function of \f r = |x| $.
\SP

\item"{(ii-c)}" There exist \f h_0 > 0 \A and, for \f 0 < h < h_0 $, a
     real-valued, measurable function \f Q_{0r}(x;\,h) \A on \f E_{R_0} \A
     such that
$$
      \sup\,\big\{\, |Q_{0r}(x;\,h)| \,/\, x \in G, \ 0 < h < h_0 \,\big\}
                                       < \infty                 \T 2.6
$$   for any compact set \f G \subset E_{R_0} $,
$$
\split
       \frac1{h}(\{Q_0((r+h)\cdot) - Q_0(r\cdot)\}\phi, \phi)
                & \le (Q_{0r}(r\cdot;\,h)\phi, \phi)  \\
                &  (\phi \in X, \ r > R_0, \ 0 < h \le h_0),
\endsplit                                                       \T 2.7
$$
     and the limit
$$
     \lim_{h\downarrow0}(Q_{0r}(r\cdot;\,h)\phi, \phi) 
          = (Q_{0r}(r\cdot)\phi, \phi)   \q (\phi \in X)        \T 2.8
$$
     exists with a real-valued, measurable, locally bounded function
     \f Q_{0r}(x) \A on \f E_{R_0} $.
\endroster
\SP
     
        (iii) There exists a positive, measurable function \f h(r) \A 
     defined on \f (R_0, \infty) \A such that
\roster
\item"{(iii-a)}" 
$$
             h(r) \le \frac{2}{r}  \qq  (r > R_0),             \T 2.9
$$
\item"{(iii-b)}" and, setting
$$
\left\{ \aligned
        & a(r) = h^{-1}(r)\sup_{|x|=r}|Q_1(x)|, \\
        & b(r) = \inf_{|x|=r} [-\,\big(Q_0(x) + h^{-1}(r)Q_{0r}(x)\big)],
\endaligned \right.                                                \T 2.10
$$
     where \f h^{-1}(r) = 1/h(r) $, we have
$$
       a(r)^2 \le b(r) \qq (r > R_0).                            \T 2.11
$$
\endroster
\BP

        In order to transform the equation (2.1) into a differential 
     equation on \f (R_0, \infty) \A \linebreak with operator-valued 
     coefficients, we give the following
\BP

        {\bf Definition 2.2.} \ (i)  For \f r > R_0 \A define a selfadjoint
     operator \f B(r) \A in \f X \A by
$$
\left\{ \aligned
             D(B(r)) = D(\Lambda_N), \\
             B(r) = -\,r^{-2}\Lambda_N,
\endaligned \right.                                               \T 2.12
$$
     where \f D(T) \A is the domain of \f T $, and \f \Lambda_N \A is the
     (selfadjoint realization of) Laplace-Beltrami operator on \f S^{N-1} $.
\SP

        (ii) For \f r > R_0 \A define a bounded operators \f C_0(r) $,
     \f C_{0r}(r;\,h) $, \f C_{0r}(r) \A and \f C_1(r) \A on \f X \A by
$$
\left\{ \aligned
              & C_0(r) = Q_0(r\cdot)\times, \\
              & C_{0r}(r;\,h) = Q_{0r}(r\cdot;\,h)\times, \\
              & C_{0r}(r) = Q_{0r}(r\cdot)\times, \\
              & C_1(r) = Q_1(r\cdot)\times.
\endaligned \right.                                                \T 2.13
$$
\BP

        {\bf Proposition 2.3.} \ {\it Let \f u \A be a solution of the
     equation {\rm (2.1)} \  and let \f v \A be as in {\rm Assumption 2.1,
     (ii)}. Let \f J = (R_0, \infty) $. Then,
\SP

        {\rm (i)} \ \f v \in C^1(J, X) $.
\SP

        {\rm (ii)} \ \f v(r) \in D((-\,\Lambda_N)^{1/2}) \A for \f r \in J $.
\SP

        {\rm (iii)} \ We have
$$
       \int_r^s \big\{ |v'(r)|^2 + |B^{1/2}(r)v(r)|^2 \big\} \, dr
                  < \infty \q (R_0 < r < s < \infty),            \T 2.14
$$
     where \f v'(r) = dv(r)/dr \A and \f B^{1/2}(r) = B(r)^{1/2} $.
\SP

        {\rm (iv)} \ \f v(r) \in D(\Lambda_N) \A for almost all
     \f r \in J $, and \f Bv \in L_2((r, s), X) \A for
     \f R_0 < r < s < \infty $.
\SP

        {\rm (v)} \ \f v'(r) \in C_{\roman ac}([r, s], X) \A for
     \f R_0 < r < s < \infty $, where \f C_{\roman ac}([r, s], X) \A is
     all \f X$-valued absolutely continuous functions on \f [r, s] $.
     There exists the weak derivative \f v''(r) \A of \f v'(r) \A for 
     \f r \in J $.
\SP

        {\rm (vi)} \ \f v'(r) \in D((-\,\Lambda_N)^{1/2}) \A for almost
     all \f r \in J $, and \f B^{1/2}v' \in L_2((r, s), X) \A for
     \f R_0 < r < s < \infty $.
\SP

        {\rm (vii)} \ \f B^{1/2}v \in C_{\roman ac}([r, s], X) \A for
     \f R_0 < r < s < \infty $, and we have
$$
\split
         \frac{d}{dr}(B^{1/2}(r)v(r), B^{1/2}(r)v(r))
              = -\frac2{r} & (B^{1/2}(r)v(r), B^{1/2}(r)v(r)) \\
              & + 2\re(B^{1/2}(r)v'(r), B^{1/2}(r)v(r))
\endsplit                                                   \T 2.15
$$
     for almost all \f r \in J $.
\SP

        {\rm (viii)} \ We have
$$
      -\,v''(r) + B(r)v(r) + C_0(r)v(r) + C_1(r)v(r) = 0        \T 2.16
$$
     in \f X \A for almost all \f r \in J $.}        
\BP

        \demo{ \ \ \ Proof} See [\YS], Proposition 1.3. \qed \enddemo
\BP

        {\bf Proposition 2.4.} \ {\it Suppose that \f Q_{0r}(x) \A satisfies
     {\rm Assumption 2.1, (ii-b)} and {\rm (ii-c)}. Let
     \f \eta \in C^1(J, X) $. Let
$$
       g(r) = \frac{d}{dr}(C_0(r)\eta(r), \eta(r))               \T 2.17
$$
     be the derivative of \f f(r) = (C_0(r)\eta(r), \eta(r)) \A in the
     sense of distributions on \f (R_0, \infty) $. Then we have
$$
       g(r) \le (C_{0r}(r)\eta(r), \eta(r))
                    + 2\re\,(C_0(r)\eta(r), \eta'(r)),         \T 2.18
$$
     where the inequality {\rm (2.18)} should be taken in the sense of
     distributions on \f (R_0, \infty) \A \linebreak again.}
\BP
     
        {\demo{ \ \ \ Proof} (I) Let \f \varphi \A be a nonnegative
     \f \Con((R_0, \infty)) \A function. Then, by definition
$$
\split
   <g, \varphi>
    & = -\, <f, \varphi'> \\
    & = -\, \lim_{h\uparrow0}
        \int_{R_0}^{\infty} f(r)\frac{\varphi(r+h)-\varphi(r)}{h}\, dr \\
    & = \lim_{h\downarrow0}\frac1{h}
         \int_{R_0}^{\infty} (f(r+h)-f(r))\varphi(r)\, dr, \\
\endsplit                                                \T 2.19
$$
     where \f < \ , \ > \A denotes the dual pair bracketing.
     
        (II) Here we have
$$
\split
    f(r+h) - f(r)
       & = (C_0(r+h)\eta(r+h), \eta(r+h)) - (C_0(r)\eta(r), \eta(r)) \\
       & = (\{C_0(r+h) - C_0(r)\}\eta(r+h), \eta(r+h)) \\
       & \ \ \ \ \ \ \ \ \ \ \ \ \
            + (C_0(r)\eta(r+h), \eta(r+h) - \eta(r)) \\
       & \ \ \ \ \ \ \ \ \ \ \ \ \ \ \ \ \ \ \
             + (\eta(r+h) - \eta(r), C_0(r)\eta(r)), \\
\endsplit                                                          \T 2.20
$$
     and hence, using (2.7) in (ii-b) of Assumption 2.1, we obtain
$$
\split
    \frac1{h}(f(r+h) - f(r))
       & \le (C_{0r}(r;\,h)\eta(r+h), \eta(r+h)) \\
       & \ \ \ \ \ \ \
            + (C_0(r)\eta(r+h), \frac1{h}(\eta(r+h) - \eta(r))) \\
       & \ \ \ \ \ \ \ \ \ \ \ \ \ \
            + (\frac1{h}(\eta(r+h) - \eta(r)), C_0(r)\eta(r)). \\
\endsplit                                                          \T 2.21
$$

        (III) It is easy to see from (2.8) in (ii-c) of Assumption 2.1 and
     (i) of Proposition 2.3 that the right-hand side of (2.21) converges to 
$$
      g_0(r) \equiv (C_{0r}(r)\eta(r), \eta(r))
              + 2\re\,(C_0(r)\eta(r), \eta(r))                  \T 2.22
$$
     boundedly on any compact interval in \f (R_0, \infty) \A as
     \f h \downarrow 0 $. Therefore, noting that \f \varphi \ge 0 $, we
     have
$$
     <g, \ \varphi> \ \le \int_{R_0}^{\infty} g_0(r)\varphi(r)\, dr
                                  = \ <g_0, \ \varphi>,        \T 2.23
$$
     which completes the proof. \qed \enddemo
$$     
     \ \ \
$$

\head  {\bf \S3. The evaluation of the functional \f M^+(v, r) $}  \endhead
\BP
  
        Let \f v = v(r\cdot) \A be as in (2.5). Then we are going to define
     the functional \f M^+(v, r) \A by
\BP

        {\bf Definition 3.1.} \ Let \f v \A be as in (ii-b) of
     Assumption 2.1. Then set
$$
      M^+(v, r) = |v'(r)|^2 - (C_0(r)v(r), v(r)) - |B^{1/2}(r)v(r)|^2
                                                                \T 3.1
$$
     for \f r > R_0 $.
\BP

        {\bf Proposition 3.2.} \ {\it Suppose that {\rm Assumption 2.1} is
     satisfied. Let \f M^+(v, r) \A be as in {\rm Definition 3.1}.
\SP

        {\rm (i)} \ Then \f M^+(v, r) \A is a real-valued, locally bounded
     function on \f J = (R_0, \infty) $. Further \f M^+(v, r) \A is right
     continuous with its left limit for \f r \in J $.
\SP

        {\rm (ii)} We have
$$
       \frac{d}{dr}M^+(v, r) \ge -\,h(r)M^+(v, r)  \qq (r > R_0),  \T 3.2
$$
     where the inequality {\rm (3.2)} \, should be taken in the sense of
     distributions on \f (R_0, \infty) $. }
\BP

        \demo{ \ \ \ Proof} (i) follows from Assumption 2.1, (ii) and
     Proposition 2.3. From Propositions 2.3 and 2.4 we see that
$$
\split
    \frac{d}{dr}M^+(v, r) \ge 2\re & (v''(r), v'(r)) \\
               & - (C_{0r}(r)v(r), v(r)) - 2\re(C_0(r)v(r), v'(r)) \\
               & + \frac2{r}(B(r)v(r), v(r)) - 2\re(B(r)v(r), v'(r))
\endsplit                                                           \T 3.3
$$
     in the sense of distributions on \f (R_0, \infty) $. Using (2.16),
     we have from (3.3)
$$
\split
    \frac{d}{dr}M^+(v, r) \ge -\,(C_{0r}(r)v(r), v(r))
                        + 2 & \re(C_1(r)v(r), v'(r)) \\
              & \ \ \ \ \ \ \ \ \ \ + \frac2{r}(B(r)v(r), v(r))   \\
              = -\, h(r)M^+(v, r) \ \ \ \ \ \ \ \ \ \ \ \ \  & \\
                 +\,h(r)\big[|v'(r)|^2 - & (C_0(r)v(r), v(r))
                                       - (B(r)v(r), v(r)) \big] \\
                -\,(C_{0r}(r)v(r), & v(r)) + 2\re(C_1(r)v(r), v'(r)) \\
                   & \ \ \ \ \ \ \ \ \ + \frac2{r}(B(r)v(r), v(r)) \\
\endsplit                                                         \T 3.4
$$
     Thus, using \f a(r) \A and \f b(r) \A defined by (2.10), and taking
     note of (2.9) in Assumption 2.1, we have
$$
\split
     \frac{d}{dr}M^+(v, r) \ge -\, h(r) & M^+(v, r) \\
          & +\, h(r)\big[\,|v'(r)|^2 -\,2a(r)|v'(r)||v(r)|
                                    +\, b(r)|v(r)|^2\,\big].
\endsplit                                                         \T 3.5
$$
     It follows from (2.11) in Assumption 2.1, that is, \f a(r)^2 \le b(r) $,
     that (3.5) implies (3.2), which completes the proof. \qed \enddemo
\BP

        {\bf Proposition 3.3.} \ {\it Suppose that {\rm Assumption 2.1} is
     satisfied. For \f R_1 > R_0 \A we have}
$$
      M^+(v, r) \ge \exp\big(-\int_{R_1}^r h(t)\,dt\,\big)M^+(v, R_1)
                                          \q (r \ge R_1).      \T 3.6
$$
\BP

        \demo{ \ \ \ Proof} It follows from Proposition 3.2 that
$$
\split
     \exp\big(& \int_{R_1}^r h(t)\,dt\,\big)\frac{d}{dr}M^+(v, r)   \\
          & + h(r)\exp\big(\int_{R_1}^r h(t)\,dt\,\big)
                           M^+(v, r) \ge 0 \q (r > R_1),
\endsplit                                                           \T 3.7
$$
     and hence
$$
    \frac{d}{dr}\bigg[\exp\big(\int_{R_1}^r h(t)\,dt\,\big)M^+(v, r)\bigg]
                                \ge 0 \q (r > R_1)               \T 3.8
$$
     in the sense of distributions on \f (R_1, \infty) $. The inequality
     (3.6) follows from (3.8) and Lemma A of Appendix. \qed \enddemo
     
$$
\ \ \ \
$$

\head {\bf \S4. The evaluation of the functional \f N(v, m, r) $} \endhead
\BP

       Using \f M^+(v, r) $, we are going to define another functional
    which will be used to evaluate \f M^+(v, r) \A in \S5.
\BP

        {\bf Definition 4.1.} \ (i) Set
$$
          N(v, m, r) = M^+(w, r) + (m(m+1) - F(r))r^{-2}|w|^2
                               \qq (w = r^mv),              \T 4.1
$$
     where \f m \A is a positive number and \f F(r) \A is a positive
     \f C^1 \A function on \f (R_0, \infty) $.
\SP

        (ii) For \f r > R_0 \A define a bounded operators \f C_R(r) \A
     on \f X \A by
$$
      C_R(r) = \re(Q(r\cdot))\times =
                (Q_0(r\cdot) + \re(Q_1(r\cdot))\times.     \T 4.2
$$
     Set
$$
          M(v, r) = |v'(r)|^2 - (C_R(r)v(r), v(r)).               \T 4.3
$$
\SP

        (iii) For \f r > R_0 \A we set
$$
        p(r) = \inf_{|x|=r}[-\,\big(2Q_0(x) + rQ_{0r}(x)\big)].   \T 4.4
$$
\BP

        {\bf Assumption 4.2.} \ The function \f F(r) \A introduced in
     Definition 4.1 satisfies the following (i) $\sim$ (iii):
\SP

     (i) There exists a positive constant \f c_0 \A such that
$$
               F^2(r) \le c_0r^4h^2(r)b(r) \qq (r > R_0),         \T 4.5
$$
     where \f b(r) \A is given in (2.10), \f F^2(r) = F(r)^2 $, and
     \f h^2(r) = h(r)^2 $.
\SP

        (ii) We have \f F(r) \to \infty \A as \f r \to \infty $.
\SP

        (iii) There exists a positive constant \f c_1 \A such that
$$
        F_r(r) \equiv \frac{d}{dr}F(r) \le c_1r^{-1}
                                 \qq (r > R_0).                  \T 4.6
$$
\BP

        {\bf Proposition 4.3.} \ (i) {\it Let \f b(r) \A be given by
     {\rm (2.10)} and assume that \f h(r) \A satisfies the inequality
     {\rm (2.9)} and that \f Q_0(x) \A is nonpositive. Then,}
$$
            r^2h^2(r)b(r)\le 2p(r) \qq  (r > R_0).               \T 4.7
$$
\SP

        (ii) {\it Assume that the inequality {\rm (2.11)} holds. Then,}
$$
   \big(r\sup_{|x|=r}|Q_1(x)|\big)^2 \le 2p(r) \qq  (r > R_0).    \T 4.8
$$
\SP

        (iii) {\it Suppose that the inequality {\rm (4,5)} holds. Then,}
$$
         r^{-2}F^2(r) \le 2c_0p(r) \qq  (r > R_0).               \T 4.9
$$
\BP

        \demo{ \ \ \ Proof} Since \f 0 < rh(r) \le 2 \A and
     \f -\,Q_0(x) \ge 0$, we have
$$
\split
      r^2h^2(r) & \big[-\,\big(Q_0(x) + h^{-1}(r)Q_{0r}(x)\big)\big] \\
                & = rh(r)\big[rh(r)\big(-\,Q_0(x)\big)
                                 + r\big(-\,Q_{0r}(x)\big)\big] \\
                & \le 2\big[2\big(-\,Q_0(x)\big)
                                  + r\big(-\,Q_{0r}(x)\big)\big] \\
                & = 2\big[-\,\big(2Q_0(x) + rQ_{0r}(x)\big)\big],
\endsplit                                                           \T 4.10
$$
     which implies (4.7). It follows from (2.11) and (4.7) that
$$
\split
     \big(r & \sup_{|x|=r}|Q_1(x)|\big)^2 \\
            & = r^2h^{2}(r)\big(h^{-1}(r)\sup_{|x|=r}|Q_1(x)|\big)^2  \\
            & \le r^2h^{2}(r)b(r)  \\
            & \le 2p(r).
\endsplit                                                      \T 4.11
$$
     From (4.5) and (4.7) we obtain
$$
    F^2(r) \le c_0r^4h^2(r)b(r) \le c_0r^2(2p(r)) = r^2(2c_0p(r))  \T 4.12
$$
     for \f r > R_0 $, which implies (4.9). \qed \enddemo
\BP

        {\bf Proposition 4.4.} \ {\it Suppose that {\rm Assumptions 2.1}
     and {\rm 4.2} hold. Then there exist \f m_0 > 0 \A and \f r_0 > R_0 \A
     such that
$$
       \frac{d}{dr}\big(r^2N(v, m, r)\big) \ge 0
                           \qq (m \ge m_0)                     \T 4.13
$$

     in the sense of distributions on \f (r_0, \infty) $.
\BP

        \demo{ \ \ \ Proof} (I) By definition \f w = r^mv \A satisfies
$$
\left\{ \aligned
        & w' = r^mv' + mr^{m-1}v = r^mv' + mr^{-1}w, \\
        & w'' = r^mv'' + 2mr^{m-1}v' + m(m-1)r^{m-2}v \\
        & \ \ \ \
             = r^mv'' + 2mr^{-1}\big(w' - mr^{-1}w\big) + m(m-1)r^{-2}w \\
        & \ \ \ \ = r^mv'' + 2mr^{-1}w' - m(m+1)r^{-2}w \\
        & \ \ \ \ = r^m(Bv + C_0v + C_1v) + 2mr^{-1}w' - m(m+1)r^{-2}w,
\endaligned \right.                                             \T 4.14
$$
     and hence we have
$$
      -\, w'' + 2mr^{-1}w'+ (B + C_0 + C_1 - m(m+1)r^{-2})w = 0.  \T 4.15
$$
\SP

        (II) Set
$$
         g(r, m) = (m(m+1) - F(r))r^{-2}.                      \T 4.16
$$
     Then, using (4.15) and Proposition 2.4, we have
$$
\split
   r^{-2}\frac{d}{dr} & \big(r^2N(v, m, r)\big) \\
     & \ge 2r^{-1}\big(|w'|^2 - (C_0w, w) - (Bw, w) + (gw, w)\big) \\
     & \ \ \ \ \ \ \ + 2\re(w'' - C_0w - Bw + gw, w')
                                + \frac{2}{r}|B^{1/2}w|^2 \\
     & \ \ \ \ \ \ \ \ \ \ \ \ \ \ \ \ \ \ \ \ \ \ \ \
                                   - (C_{0r}w, w) + (g_rw, w) \\
     & = 2r^{-1}\big(|w'|^2 - (C_0w, w) - (Bw, w) + (gw, w)\big) \\
     & \ \ \ \ \ \ \ + 2\re(2mr^{-1}w' + C_1w - m(m+1)r^{-2}w + gw, w') \\
     & \ \ \ \ \ \ \ \ \ \ \ \ \ \ \ \ \ \ \
                 + \frac{2}{r}|B^{1/2}w|^2- (C_{0r}w, w) + (g_rw, w) \\
     & = 2(1 + 2m)r^{-1}|w'|^2 + (2r^{-1}g + g_r)|w|^2 \\
     & \ \ \ \ \ \ \ + r^{-1}([-\,2C_0 - rC_{0r}]w, w) \\
     & \ \ \ \ \ \ \ \ \ \ \ \ \ \ \ \
                         + 2\re(C_1w - m(m+1)r^{-2}w + gw, w'),
\endsplit                                                           \T 4.17
$$
     where \f g_r = dg/dr $. Note that
$$
\left\{ \aligned
     & 2r^{-1}g(r, m) + g_r(r, m) = -\,F_r(r)r^{-2},  \\
     & g(r, m) - m(m+1)r^{-2} = -\,F(r)r^{-2}.
\endaligned \right.                                                \T 4.18
$$
     Then the above inequality (4.17) can be rewritten as
$$
\split
      r^{-2}\frac{d}{dr} & \big(r^2N(v, m, r)\big) \\
          & \ge 2(1 + 2m)r^{-1}|w'|^2
                  + r^{-1}(p(r) - r^{-1}F_r(r))|w|^2 \\
          & \ \ \ \ \ \ \ \ \ \ \ \ \ \ \ \ \ \ \ \ \ \ \ \ \ \ \ \
                   - 2\re((r^{-2}F(r) - C_1)w, w'),
\endsplit                                                        \T 4.19
$$
     where \f p(r) \A is as in (4.4).
\SP

        (III) It follows from Assumption 4.2, (iii) and Proposition 4.3,
     (ii), (iii) that
$$
\left\{ \aligned
      & p(r) - r^{-1}F_r(r) \ge p(r) - c_1r^{-2} = r^{-2}(r^2p(r) - c_1), \\
      & |r^{-2}F(r) - Q_1(x)| \le r^{-1}\big(2c_0p(r)\big)^{1/2} 
                          +  r^{-1}\big(2p(r)\big)^{1/2}     \\
      & \ \ \ \ \ \ \ \ \ \ \ \ \ \ \ \ \ \ \ \ \
            = c_2r^{-1}p^{1/2}(r)  \q (c_2 = \sqrt{2c_0} + \sqrt{2}).
\endaligned \right.                                              \T 4.20
$$
     Further, from (iii) of Proposition 4.3 and (ii) of Assumption 4.2 we
     see that
$$
     r^2p(r) \ge (2c_0)^{-1}F^2(r) \to \infty \q (r \to \infty), \T 4.21
$$
     and hence there exists \f r_0 > R_0 \A such that
$$
     r^{-2}(r^2p(r) - c_1) = p(r)(1 - \frac{c_1}{r^2p(r)})
                                  \ge 2^{-1}p(r)              \T 4.22
$$
     for \f r \ge r_0 $. Thus, we obtain from (5.19)
$$
\split
     r^{-2}\frac{d}{dr} & \big(r^2N(v, m, r)\big) \\
          & \ge 2(1 + 2m)r^{-1}|w'|^2 + 2^{-1}r^{-1}p(r)|w|^2 \\
          & \ \ \ \ \ \ \ \ \ \ \ \ \ \ \ \ \ \ \ \ \ \ \ \ \ \ \ \
                           - 2c_2r^{-1}p^{1/2}(r)|w'||w| \\
          & = r^{-1}\big[2(1 + 2m)|w'|^2 + 2^{-1}p(r)|w|^2  \\
          & \ \ \ \ \ \ \ \ \ \ \ \ \ \ \ \ \ \ \ \ \ \ \ \ \ \ \ \
                               - 2c_2p^{1/2}(r)|w'||w|\big]
\endsplit                                                        \T 4.23
$$
     for \f r \ge r_0 $. Therefore, there exists a sufficiently large
     \f m_0 > 0 \A such that
$$
       r^2\frac{d}{dr}\big(r^2N(v, m, r)\big) \ge 0                 \T 4.24
$$
     for \f r \ge r_0 \A and \f m \ge m_0 $, which completes the proof.
     \qed \enddemo
     
\newpage

%$$
%\ \ \
%$$

\head      {\bf \S5. Uniqueness Theorem}        \endhead
\BP

        We are going to prove our main theorem (Theorem 5.10) which shows,
     under Assumptions 2.1 and 4.1, and some additional conditions
     (Assumptions 5.5 and 5.8), that the solution \f u \A has compact
     support if \f u \A satisfies
$$
       \varliminf_{r\to\infty}
           \int_{|x|=r} \big\{ \big|\frac{\pa u}{\pa r}\big|^2
                - \re\,(q(x))|u|^2 \big\}\, dS = 0.        \T 5.1
$$
\BP

        {\bf Proposition 5.1.} \ {\it Suppose that {\rm Assumptions 2.1} and
     {\rm 4.2} hold. Suppose that the support of \f u \A is unbounded. Let
     \f r_0 \A and \f m_0 \A be as in} Proposition 4.4. {\it Then
     there exist \f m_1 \ge m_0 \A and \f r_1 \ge r_0 \A such that}
$$
        N(v, m_1, r) > 0   \qq (r \ge r_1).                    \T 5.2
$$
\BP

        \demo{ \ \ \ Proof} Since the support of \f u \A is assumed to be
     unbounded, there exists \f r_1 \ge r_0 \A such that \f |v(r_1)| > 0 $. 
     Since
$$
\split
    r_1^{-2m} & N(v, m, r_1)  \\
      & = r_1^{-2m}\big\{\,|w'(r_1)|^2 - (C_0(r_1)w(r_1), w(r_1))  \\
      & \ \ \ \ \ \ \ \ \ \ \
      - |B^{1/2}(r_1)w(r_1)|^2 + \big(m(m+1) - F(r_1)\big)|w(r_1)|^2\big\} \\
      & \ge -\,(C_0(r_1)v(r_1), v(r_1))  \\
      & \ \ \ \ \ \ \ \ \ \ - |B^{1/2}(r_1)v(r_1)|^2
                   + \big(m(m+1) - F(r_1)\big)|v(r_1)|^2,
\endsplit                                                        \T 5.3
$$
        we can choose a sufficiently large \f m_1 \ge m_0 \A so that
$$
        r_1^{-2m_1}N(v, m_1, r_1) > 0, \ \  \roman{or} \ \
                                  r_1^2 N(v, m_1, r_1) > 0.      \T 5.4
$$
     Note that, by (ii)-2 of Assumption 2.1, \f N(r, m, v) \A is
     right-continuous. Then the inequality (5.4) is combined with (4.13) 
     and Lemma A in Appendix to see that \f r^2N(r, m_1, v) > 0 \A on
     \f [r_1, \infty) $, which completes the proof. \qed \enddemo
\BP

        {\bf Definition 5.2.} \ Suppose that Assumptions 2.1 and 4.2 hold.
     Suppose that the support of \f u \A is unbounded. Let \f F(r) \A and
     \f m_1 \A be given in Definition 4.1 and Proposition 5.1, respectively.
     Then we introduce the following two alternative cases:
\SP

        Case I\,: There exists an infinite sequence \f \{\,r_{\ell}'\,\} \A
     such that \f R_0 < r_{\ell}' $, \f r_{\ell}' \to \infty \A as
     \f \ell \to \infty $, and
$$
       2\re(v'(r_{\ell}'), v(r_{\ell}'))
          \le (2m_1r_{\ell}')^{-1}F(r_{\ell}')|v(r_{\ell}')|^2   \T 5.5
$$
     for all \f \ell = 1, 2, \cdots $.
\SP

        Case II\,: There exists \f r_2 > r_1 \A such that
$$        
        2\re(v'(r), v(r)) > (2m_1r)^{-1}F(r)|v(r)|^2.
                        \qq (r \ge r_2),                         \T 5.6
$$
     where \f r_1 \A is as in Proposition 5.1.
\BP

        {\bf Proposition 5.3.} \ {\it Suppose that {\rm Assumptions 2.1} and
     {\rm 4.2} hold. Suppose that the support of \f u \A is unbounded.
     Suppose that {\rm Case I} in {\rm Definition 5.2} holds. Then there
     exists an infinite sequence \f \{\,r_{\ell}''\,\} \A such that
     \f R_0 < r_{\ell}'' $, \A r_{\ell}'' \to \infty \A as
     \f \ell \to \infty $, and}
$$
       M^+(v, r_{\ell}'') > 0 \qq (\ell = 1, 2, \cdots).          \T 5.7
$$
\BP

        \demo{ \ \ \ Proof} Let \f \{\,r_{\ell}'\,\} \A be as in Case I
     of Definition 5.2. Let \f w = r^{m_1}v $, where \f m_1 \A is as in
     Proposition 5.1. Then we have for\f r = r_{\ell}' \A
$$
\split
     r^{-2m_1}|w'|^2 & = |v' + m_1r^{-1}v |^2 \\
                     & = |v'|^2 + 2m_1r^{-1}\re(v', v) + m_1^2r^{-2}|v|^2 \\
                     & \le |v'|^2 + m_1r^{-1}(2m_1r)^{-1}F(r)|v|^2
                                                 + m_1^2r^{-2}|v|^2 \\
                     & = |v'|^2 + \big(2^{-1}F(r) + m_1^2\big)r^{-2}|v|^2.
\endsplit                                                         \T 5.8
$$
     Let \f r_1 \A be as in Proposition 5.1. For \f r = r_{\ell}' \A such
     that \f r_{\ell}' \ge r_1 $, it follows that
$$
\split
    0 & < N(v, m_1, r) = M^+(w, r) + (m_1(m_1+1) - F(r))r^{-2}|w|^2 \\
      & \le r^{2m_1}
           \big\{|v'|^2 + \big(2^{-1}F(r) + m_1^2\big)r^{-2}|v|^2\big\} \\
      & \ \ \ \ \ \ \ \ \ \ \ \ \ \ \ \ \
             - r^{2m_1}\big\{(C_0v, v) + (Bv, v)\big\} \\
      & \ \ \ \ \ \ \ \ \ \ \ \ \ \ \ \ \
             + r^{2m_1}(m_1(m_1+1) - F(r))r^{-2}|v|^2 \\
      & = r^{2m_1}
          \big\{M^+(v, r) + (m_1(2m_1+1) - 2^{-1}F(r))r^{-2}|v|^2\big\}.
\endsplit                                                           \T 5.9
$$
     Since \f F(r) \to \infty \A as \f r \to \infty $, there exists
     a positive integer \f \ell_0 \A such that
$$
       m_1(2m_1+1) - 2^{-1}F(r_{\ell}')) < 0
                                \qq (\ell \ge \ell_0).             \T 5.10
$$
     Therefore we have only to define \f r_{\ell}'' \A by
$$
       r_{\ell}'' = r_{\ell_0 + \ell}' \q (\ell = 1, 2, \cdots),    \T 5.11
$$
     which completes the proof. \qed \enddemo
\BP

        {\bf Proposition 5.4.} \ {\it Suppose that {\rm Assumptions 2.1}
     and {\rm 4.2} hold. Suppose that the support of \f u \A is unbounded.
     Suppose that {\rm Case II} in {\rm Definition 5.2} holds. Suppose,
     in addition, that
$$
          \re\, Q(x) \le 0 \qq (x \in E_{R_0}).                   \T 5.12
$$
     Then there exist \f r_3 > R_0 \A and a positive constant \f c_2 \A
     such that
$$
             M(v, r) \ge c_2 \qq (r \ge r_3),                    \T 5.13
$$
     where \f M(v, r) \A is given by {\rm (4.3)}.}
\BP

        \demo{ \ \ \ Proof} Since \f F(r) \to \infty \A as \f r \to \infty $,
     there exists \f r_4 > R_0 \A such that
$$
         \frac{F(r)}{2m_1} \ge 2 \qq (r \ge r_4).                \T 5.14
$$
     Then it follows from (5.6) that
$$
         \frac{d}{dr}|v(r)|^2 \ge 2r^{-1}|v(r)|^2\q (r \ge r_4).  \T 5.15
$$
     Let \f r_3 \A be such that \f r_3 \ge r_4 \A and \f |v(r_3)| > 0 $.
     Then, since
$$
    \frac{d}{dr}\big(r^{-2}|v(r)|^2\big) =
          r^{-2}\big(\frac{d}{dr}|v(r)|^2 - 2r^{-1}|v(r)|^2\big) \ge 0
                                           \q (r \ge r_4),        \T 5.16
$$
     we have
$$
        r^{-2}|v(r)|^2 \ge r_3^{-2}|v(r_3)|^2 > 0 \q
                                     (r \ge r_3).                \T 5.17
$$
     Also, using (5.6) and (5.14) again, we see that
$$
       2r^{-1}|v(r)|^2 \le (2m_1r)^{-1}F(r)|v(r)|^2
                   \le 2|v(r)||v'(r)|,                     \T 5.18
$$
     or
$$
        r^{-1}|v(r)| \le |v'(r)|                              \T 5.19
$$
     for \f r \ge r_4 $. Thus, it follows from (5.17) and (5.19) that
$$
      |v'(r)|^2 \ge r^{-2}|v(r)|^2 \ge r_2^{-2}|v(r_3)|^2 > 0 \q \T 5.20
$$
     for \f r \ge r_3 $, which is combined with (5.12) to obtain (5.13).
     \qed \enddemo
\BP

        {\bf Assumption 5.5.} \ (i) Let \f h(r) \A be as above. Then
     \f h \in L_1((R_0, \infty)) $.
\SP

        (ii) There exists a constant \f \beta \in (0, 1) \A such that
$$
        0 \ge \beta Q_0(x) \ge \re\,(Q(x)) \q  (x \in E_{R_0}).   \T 5.21
$$
\BP

        {\bf Theorem 5.6.} \ {\it Suppose that {\rm Assumptions 2.1, 4.2}
     and {\rm 5.5} hold. Suppose that the support of \f u \A is unbounded.
     Then there exist a positive constant \f c_3 \A  and \f R_2 > R_0 \A
     such that}
$$
          M(v, r) \ge c_3 \qq (r \ge R_2).                   \T 5.22
$$
\BP

        \demo{ \ \ \ Proof} Note that all the assumptions that are
     necessary for the conclusions of Propositions 3.2, 3.3, 4.3, 4.4, 5.1
     5.3 and 5.4 are satisfied. Suppose that Case I of Definition 5.2 is
     satisfied. Then, by Proposition 5.3 there exists \f R_2' > R_0 \A such
     that \f M^+(v, R_2') > 0 $. Therefore, setting \f R_1 = R_2' \A
     in Proposition 3.3, we have for \f r \ge R_2' $,
$$
\split
    M^+(v, r) & \ge \exp\big(-\int_{R_3'}^r h(t)\,dt\,\big)M^+(v, R_2') \\
              & \ge \exp\big(-\int_{R_3'}^{\infty}
                    h(t)\,dt\,\big)M^+(v, R_2').  \\
\endsplit                                                         \T 5.23
$$
     Since we have from (5.21)
$$
\split
    M(v, r) & = |v'(r)|^2 - (C_R(r)v(r), v(r)) \\
            & \ge |v'(r)|^2 - \beta(C_0(r)v(r), v(r)) \\
            & \ge \beta|v'(r)|^2 - \beta(C_0(r)v(r), v(r))
                                 - \beta|B^{1/2}(r)v(r)|^2 \\
            & = \beta M^+(v, r),
\endsplit                                                         \T 5.24
$$
     it follows from (5.23) that
$$
         M(v, r) \ge c_3' \qq (r \ge R_2')                         \T 5.25
$$
     with
$$
       c_3' = \beta \exp\big(-\int_{R_3'}^{\infty}
                                    h(t)\,dt\,\big)M^+(v, R_2').  \T 5.26
$$
     Suppose that Case II of Definition 5.2 is satisfied. Then from
     Proposition 5.4 we have
$$
             M(v, r) \ge c_2 \qq (r \ge r_3),                      \T 5.27
$$
     where \f c_2 \A and \f r_3 \A are as in Proposition 5.4. Now set
$$
\left\{ \aligned
        & c_3 = \min\,\{c_3', \, c_2\}, \\
        & R_2 = \max\,\{R_2', \, r_3\}.
\endaligned \right.                                                \T 5.28
$$
     Then (5.22) follows, which completes the proof. \qed \enddemo
\BP

        {\bf Corollary 5.7.} \ {\it Suppose that {\rm Assumptions 2.1, 4.2}
     and {\rm 5.8} hold. Suppose that
$$
             \varliminf_{r\to\infty} M(v, r) = 0.                 \T 5.29
$$
     Then \f u \A has compact support.}
\BP

        In order to show our main theorem (Theorem 5.10) we need one more
     assumption.
\BP        
        
        {\bf Assumption 5.8.} \ We have
$$
      \lim_{r\to\infty}\big(r^2\inf_{|x|=r}\re\,(-q(x))\big) = \infty.
                                                               \T 5.30
$$
\BP

        Before we state and prove Theorem 5.10, we are going to unify 
     Assumptions 2.1, 4.2, 5.5 and 5.8 in more organized form:
\BP

        {\bf Assumption 5.9.} \ (1) Let \f N \A be an integer such that
     \f N \ge 2 $. Let \f u \in H^2(E_{R_0})_{\roman loc} $, \f R_0 > 0 $,
     be a solution of the homogeneous \Sch\ equation (2.1), where 
     \f E_{R_0} \A is given by (2.2) (with \f R = R_0 $). Here \f q(x) \A 
     is a complex-valued, measurable, locally bounded function on 
     \f E_{R_0} \A which satisfies (5.30).
\SP

        (2) Set
$$
           Q(x) = q(x) + \frac{(N-1)(N-3)}{4r^2}.             \T 5.31
$$
\SP

\roster
\item"{(2-a)}" Then \f Q(x) \A is decomposed as
$$
         Q(x) = Q_0(x) + Q_1(x),                             \T 5.32
$$
     where \f Q_0(x) \A is a non-positive, measurable, locally bounded
     function on \f E_{R_0} \A \linebreak and \f Q_1(x) \A is a 
     complex-valued, measurable, locally bounded function on \f E_{R_0} \A 
     such that
$$
        0 \ge \beta Q_0(x) \ge \re\,(Q(x)) \q  (x \in E_{R_0})   \T 5.33
$$
     with a constant \f \beta \in (0, 1) $. 
\SP

\item"{(2-b)}" For any \f \phi \in X = L_2(S^{N-1}) $, 
     \f (Q_0(r\cdot)\phi, \phi) \A has the right limit for all \f r > R_0 \A
     as a function of \f r = |x| $, where \f ( \ , \ ) \A is the inner 
     product of \f X $.
\SP

\item"{(2-c)}" There exist \f h_0 > 0 \A and, for \f 0 < h < h_0 $, a
     real-valued, measurable function \f Q_{0r}(x;\,h) \A on \f E_{R_0} \A
     such that (2.6), (2.7) and (2.8) hold.
\SP
     
\item"{(2-d)}" There exists \f h(r) \in L_1((R_0, \infty)) \A such that
$$
          0 < h(r) \le \frac{2}{r}  \qq  (r > R_0),             \T 5.34
$$
     and, setting 
$$
\left\{ \aligned
        & a(r) = h^{-1}(r)\sup_{|x|=r}|Q_1(x)|, \\
        & b(r) = \inf_{|x|=r} [-\,\big(Q_0(x) + h^{-1}(r)Q_{0r}(x)\big)],\\
        & \ \ \ \ \ \ \ \ \ \ \ \qq \qq (h^{-1}(r) = 1/h(r)),
\endaligned \right.                                                \T 5.35
$$
     we have
$$
       a(r)^2 \le b(r) \qq (r > R_0).                            \T 5.36
$$
\endroster
\SP

        (3) The function \f F(r) \A introduced in Definition 4.1 satisfies
      Assumption 4.2.
\BP
        
        {\bf Theorem 5.10.} \ {\it Suppose that {\rm Assumptions 5.9} hold. 
     Suppose that the solution \f u \A satisfies {\rm (5.1)}. Then 
     \f u \A has compact support.}
\BP

        \demo{ \ \ \ Proof} Note that
$$
\split
   0 \le |v'(r)|^2 & - (C_R(r)v(r), v(r)) \\
                   & = |\big(r^{(N-1)/2}u(r\cdot)\big)'|^2
                        - r^{N-1}(\re\,(Q(r\cdot))u(r\cdot), u(r\cdot)) \\
                   & = r^{N-1}|\pa_ru(r\cdot)
                         + 2^{-1}(N-1)r^{-1}u(r\cdot)|^2  \\
                   & \ \ \ \ \ \ \ \ \ \ \ \ \ \ \
                       - r^{N-1}(\re\,(Q(r\cdot))u(r\cdot), u(r\cdot)) \\
                   & \le 2r^{N-1}|\pa_ru(r\cdot)|^2
                               + 2^{-1}(N-1)^2r^{N-3}|u(r\cdot)|^2 \\
                   & \ \ \ \ \ \ \ \ \ \ \ \ \ \ \
                            - r^{N-1}(\big\{\re\,(q(r\cdot))
                                 + \frac{(N-1)(N-3)}{4r^2}\big\}
                                        u(r\cdot), u(r\cdot)), \\
                   & \le 2r^{N-1}|\pa_ru(r\cdot)|^2
                          + \frac{N^2-1}4\,r^{N-3}|u(r\cdot)|^2  \\
                   & \ \ \ \ \ \ \ \ \ \ \ \ \ \ \ \ \ \ \ \ \ \ \
                        - r^{N-1}(\re\,(q(r\cdot))u(r\cdot), u(r\cdot)),
\endsplit                                                          \T 5.37
$$
     where \f \pa_r = \pa/\pa r \A and we have used
$$
\left\{ \aligned
           & Q(x) = q(x) + \frac{(N-1)(N-3)}{4r^2},     \\
           & \frac{(N-1)^2}2 - \frac{(N-1)(N-3)}{4} = \frac{N^2-1}4.
\endaligned \right.                                               \T 5.38
$$
     Therefore we have
$$
\split
    0 \le |v'(r)|^2 & - (C_R(r)v(r), v(r)) \\
       \le & \ 2r^{N-1}\big\{|\pa_ru(r\cdot)|^2
            - (\re\,(q(r\cdot))u(r\cdot), u(r\cdot))\big\} \\
           & \ \ \ - r^{N-3}(\big[-r^2\re\,(q(r\cdot)) - 4^{-1}(N^2-1)\big]
                                       u(r\cdot), u(r\cdot)). \\
\endsplit                                                         \T 5.39
$$
     It follows from (1) of Assumption 5.9 ((5.30)) that there exists
     \f R_3 > R_0 \A such that
$$
      -\,r^2\re\,(q(x)) - 4^{-1}(N^2-1) > 0 \qq (|x| = r, \ r \ge R_3),
                                                                \T 5.40
$$
     and hence, for \f |x| \ge R_3 $,
$$
\split
    0 \le |v'(r)|^2 & - (C_R(r)v(r), v(r)) \\
        & \le 2\int_{|x|=r} \big\{ \big|\frac{\pa u}{\pa r}\big|^2
                 - \re\,(q(x))|u|^2 \big\}\, dS  , \\
\endsplit                                                           \T 5.41
$$
     which, together with (5.1), implies (5.29). Thus Corollary 5.7 can be
     applied to see that \f u \A has compact support, which completes
     the proof. \qed \enddemo

$$
\ \ \
$$

\newpage

\head     {\bf \S6. Examples}                    \endhead
\BP
 
        In this section we are going to give some applications of
     Theorem 5.10.
\BP

        {\bf Example 6.1.} \ Let \f \overline{R} > 0 \A and let
     \f u \in H^2(E_{\overline{R}})_{\roman loc} \A be a solution of the
     equation
$$
      (-\,\Delta + V_{\ell}(x) + V_s(x) - \lambda(x))u = 0 \q
                               (x \in E_{\overline{R}})         \T 6.1
$$
     Here \f \lambda(x) \A is a real-valued, measurable, locally bounded
     function on \f E_{\overline{R}} \A \, satisfying the following (i) and
     (ii):
\roster
\item"{(i)}" There exists \f m_0 > 0 \A such that
$$
           \lambda(x) \ge m_0 \qq      (x \in E_{\overline{R}}). \T 6.2
$$
\item"{(ii)}" For any \f \phi \in X \A the function
$$
          f_{\phi}(r) = (\lambda(r\cdot)\phi, \phi)               \T 6.3
$$
     is a right continuous, nondecreasing function on
     \f (\overline{R}, \infty) $.
\endroster
\F   The functions \f V_{\ell}(x) \A and \f V_s(x) \A are real-valued and
     complex-valued functions, respectively, satisfying the following
     (iii) and (iv):
\roster
\item"{(iii)}" The long-range potential \f V_{\ell}(x) \A is assumed to
     be \f C^1 \A function on \f E_{\overline{R}} \A such that
$$
\left\{ \aligned
    & \lim_{r\to\infty}\,\sup_{|x|=r} |V_{\ell}(x)| = 0, \\
    & \sup_{r>\overline{R},\,|x|=r}
           \big\{\vert x \vert^{1+\epsilon} \big|
             \frac{\pa V_{\ell}}{\pa \vert x \vert}\big|\big\} < \infty
\endaligned \right.                                            \T 6.4
$$
     with \f \epsilon \in (0, 2) $.
\item"{(iv)}" The short-range potential \f V_s(x) \A  is assumed to
     be measurable such that
$$
     \sup_{r>\overline{R},\,|x|=r}
             \big\{r^{1+\epsilon}|V_{s}(x)|\big\} < \infty,     \T 6.5
$$
     where \f \epsilon \A is as above.
\endroster
\F   Suppose, in addition, that
$$
     \varliminf_{r\to\infty}
          \int_{|x|=r} \big\{ \big|\frac{\pa u}{\pa r}\big|^2
                     + \lambda(x)|u|^2 \big\}\, dS = 0.          \T 6.6
$$
     Then \f u \A is identically zero in \f E_{\overline{R}} $. In fact, set
$$
\left\{ \aligned
     & Q_0(x) = -\,\lambda(x) + V_{\ell}(x), \\
     & Q_{0r}(r\omega;\,h)
             = h^{-1}\int_r^{r+h} \frac{\pa V_{\ell}}{\pa r}(s\omega)\,ds
                  \q (\omega \in S^{N-1}, \, h > 0), \\
     & Q_{0r}(x) = \frac{\pa V_{\ell}}{\pa r}, \\
     & Q_1(x) = V_s(x) + \frac{(N-1)(N-3)}{4r^2}.
\endaligned \right.                                       \T 6.7
$$
     Then, since
$$
       \re\,(-q(x)) = \lambda(x) - V_{\ell}(x) - \re\,(V_s(x))
                  \ge m_0 - V_{\ell}(x) - \re\,(V_s(x))
                                                      \T 6.8
$$
     (1) of Assumption 5.9 is satisfied for sufficiently large \f r $.
     For \f \phi \in X $, we have
$$
\split
     \frac1{h} & (\big[Q_0((r+h)\cdot) - Q_0(r\cdot)\big]\phi, \phi)  \\
      & \ \ \ \ \ \ \ \ = - \frac1{h}
          (\big[\lambda((r+h)\cdot) - \lambda(r\cdot)\big]\phi, \phi)
           + \frac1{h}
             (\big[V_{\ell}((r+h)\cdot) - V_{\ell}(r\cdot)\big]\phi, \phi) \\
      & \ \ \ \ \ \ \ \
              \le (Q_{0r}(r\cdot;\,h)\phi, \phi)              \\
      & \ \ \ \ \ \ \ \  \to (Q_{0r}(r\cdot)\phi, \phi)
\endsplit                                                          \T 6.9
$$
     as \f h \to 0 \A with \f h > 0 $. Thus (2-c) of Assumption 5.9 is
     satisfied. Set
$$
       h(r) = r^{-1-\epsilon/2}  \q (r > \overline{R}).         \T 6.10
$$
     Then \f h(r) \in L_1((\overline{R}, \infty)) \A and the inequality
     (5.34) is satisfied for sufficiently large \f r $. Also we have
     \f a(r) \to 0 \A as \f r \to \infty \A and \f b(r) \ge m_0/2 \A for
     sufficiently large \f r $, and hence (2-d) of Assumption 5.9 is now 
     satisfied. Noting that
$$
      \beta Q_0(x) - \re\,(Q(x)) =
           (1 - \beta)\lambda(x) + (\beta - 1)V_{\ell}(x)
                                      - \re\,(V_s(x)),      \T 6.11
$$
     and that \f \lambda(x) \ge m_0 \A ((6.2)), we see that (2-a) of 
     Assumption 5.9 holds for sufficiently large \f r \A with any
     \f \beta \in (0, 1) $. The condition (2-b) of Assumption 5.9 is 
     verified by (ii) of Example 6.1 and the smoothness of \f V_{\ell}(x) $.
     Define \f F(r) \A by \f F(r) = \log r $. Obviously (ii) and (iii)
     of Assumption 4.2 are satisfied by definition. Since
$$
          r^4h^2(r)b(r) = r^{2-\epsilon}(\lambda(x) + o(1))     \T 6.12
$$
     as \f r \to \infty $, (4.5) in Assumption 4.2 holds for sufficiently
     large \f r $. Therefore, by setting \f R_0 \A sufficiently large,
     all the conditions of Assumption 5.9 are satisfied, which implies that
     the solution \f u \A has compact support in \f E_{\overline{R}} $.
     Therefore it follows from the unique continuation theorem that \f u \A
     is identically zero in \f E_{\overline{R}} $.
     
        We remark here that, if \f \lambda(x) \A is assumed to be bounded
     from above, too, then the condition (6.6) is equivalent to
$$
      \varliminf_{r\to\infty}
           \int_{|x|=r} \big\{ \big|\frac{\pa u}{\pa r}\big|^2
                    + |u|^2 \big\}\, dS = 0.                    \T 6.13
$$
     Another remark is that, if \f V_s(x) \A is real-valued, then
     the condition (6.6) is implied by the generalized radiation condition
$$
       r^{\delta-1}\big(\frac{\pa u}{\pa r} - i\sqrt{\lambda(x)}u\big)
                                            \in L_2(E_R)       \T 6.14
$$
     with \f \delta > 1/2 \A and \f R > \overline{R} $.
\BP

        {\bf Example 6.2.} \ Let \f \overline{R} > 0 \A and let
     \f u \in H^2(E_{\overline{R}})_{\roman loc} \A be a solution of the
     equation
$$
       (-\,\frac1{\mu(x)}\Delta - \lambda)u = 0 \q
                                 (x \in E_{\overline{R}})    \T 6.15
$$
     Here \f \lambda > 0 \A and the real-valued function \f \mu(x) \A on
     \f E_{\overline{R}} \A is decomposed as
$$
         \mu(x) = \mu_0(x) + \mu_{\ell}(x) + \mu_s(x)
                               \q (x \in E_{\overline{R}}),    \T 6.16
$$
     where \f \mu_0(x) $, \f \mu_{\ell}(x) \A and \f \mu_s(x) \A satisfy
     the following (i)$\sim$(iv):
\roster
\item"{(i)}" \f \mu_0(x) \A is real-valued and measurable and there exists
     \f \widetilde{m_0} > 0 \A such that
$$
        \mu_0(x) \ge \widetilde{m_0} \qq
                                  (x \in E_{\overline{R}}).    \T 6.17
$$
\item"{(ii)}" For any \f \phi \in X \A the function
$$
           g_{\phi}(r) = (\mu_0(r\cdot)\phi, \phi)               \T 6.18
$$
     is a right continuous, nondecreasing function on
     \f (\overline{R}, \infty) $.
\item"{(iii)}" The real-valued function \f \mu_{\ell} \A satisfies
     Example 6.1, (iii) with \f V_{\ell}(x) \A replaced by
     \f \mu_{\ell} $.
\item"{(iv)}" The complex-valued function \f \mu_s(x) \A  satisfies
     Example 6.1, (iv) with \f V_s(x) \A re\-placed by \f \mu_s $.
\endroster
\F   Set
$$
\left\{ \aligned
         & \lambda(x) = \lambda\mu_0(x), \\
         & V_{\ell}(x) = \lambda\mu_{\ell}(x), \\
         & V_s(x) = \lambda\mu_s(x)
\endaligned \right.                                              \T 6.19
$$
     Then \f u \A satisfy the equation (6.1) in Example 6.1, where
     \f \lambda(x) $, \f V_{\ell}(x) $ and \f V_s(x) \A \linebreak satisfy
     (i)$\sim$(iv) in Example 6.1. Thus the condition
$$
     \varliminf_{r\to\infty}
         \int_{|x|=r} \big\{ \big|\frac{\pa u}{\pa r}\big|^2
               + \lambda\mu_0(x)|u|^2 \big\}\, dS = 0        \T 6.20
$$
     implies that \f u \A is identically zero.
     
$$
\ \ \ \ \
$$

\head  {\bf \S7 Reduced Wave operator in layered media}   \endhead
\BP

        In [\JSb] we considered the reduced wave operator
$$
            H = - \frac1{\mu(x)}\Delta \qq \roman{in \ } 
                                  \CH = L_2(\RN, \mu(x)dx),         \T 7.1
$$
     where \f \mu(x) \A is a real-valued function such that 
$$
       0 < \inf_x \mu(x) \le \sup_x \mu(x) < \infty.                \T 7.2
$$
     By defining the domain \f D(H) \A of \f H \A by \f D(H) = H^2(\RN) $,
     where \f H^2(\RN) \A is the second order Sobolev space on \f \RN $,
     \f H \A becomes a self-adjoint operator on \f \CH $. In this section
     we shall show that the nonexistence of the eigenvalues of \f H \A can
     be proved in some cases discussed in [\JSb] by using the result of \S6
     (Example 6.2). Suppose that \f \mu(x) \A has the decomposition (6.16)
     with a positive function \f \mu_0 $, a long-range perturbation
     \f \mu_{\ell} \A and a short-range perturbation \f \mu_s $. The 
     functions \f \mu_0 $, \f \mu_{\ell} \A and \f \mu_s \A are assumed to
     satisfy (i) $\sim$ (iv) of Example 6.2. In [\JSb], for the sake of
     simplicity, we assumed that only one of a long-range perturbation or
     short-range perturbation appeared with the main term \f \mu_0(x) $, 
     but we can easily modify the arguments in [\JSb] so that we can treat 
     \f \mu(x) \A of the form (6.16). Let \f K_- \A be a nonpositive integer
     or \f K_- = -\infty \A and let \f K_+ \A be a nonnegative integer or
     \f K_+ = \infty $. Let \f K \A be a set of integers given by
$$
             K = \{ k /  K_- \le k \le K_+ \}.            \T 7.3
$$
     Let \f \{ \Omega_k \}_{k \in K} \A be a sequence of open sets of
     \f \RN \A such that
$$
\left\{ \aligned
          & \Omega_k \cap \Omega_{\ell} = \emptyset \q (k \ne \ell),   \\
          & \bigcup_{k\in K} \overline{\Omega_k} = \RN,
\endaligned \right.                                              \T 7.4
$$
     where \f \overline{A} \A is the closure of \f A $. Further
     we assume that the boundary \f \pa\Omega_k \A of \f \Omega_k \A has the 
     form 
$$
             \pa\Omega_k = S_k^{(-)} \cup S_k^{(+)},               \T 7.5
$$ 
     where \f S_k^{(-)} \cap S_k^{(+)} = \emptyset $, and each of 
     \f S_k^{(-)}\A and \f S_k^{(+)} \A is a continuous surface which is a 
     finite union of smooth surfaces. We also assume that
$$
\left\{ \aligned
          & S_k^{(+)} = S_{k+1}^{(-)} \qq \ \ (k \in K), \\
          & S_{K_+}^{(+)} = S_{K_++1}^{(-)} = \emptyset
                                      \q (\roman{if \ } K_+ \ne \infty), \\
          & S_{K_-}^{(-)} = \emptyset 
                    \qq \ \ \ \ \ \ \  (\roman{if \ } K_- \ne -\infty). 
\endaligned \right.                                                 \T 7.6
$$
     Now the function \f \mu_0(x) \A is assumed to be a simple function
     which takes a constant value \f \nu_k \A on each \f \Omega_k \A such
     that \f \{ \nu_k \}_{k\in K}\A ia a bounded, positive sequence. We 
     assume that the origin \f 0 \A of the coordinates is in \f \Omega_0 $,
     and \f \mu_0(x) \A satisfies the condition
$$
        (\nu_{k+1} - \nu_k)(n^{(k)}(x)\cdot x) \ge 0 
                  \q (x \in S_k^{(+)} = S_{k+1}^{(-)}, \ k \in K),   \T 7.7
$$
     where \f n^{(k)}(x) \A is the unit outward normal of \f \Omega_k \A at 
     \f x \in \pa\Omega_k \A and \f n^{(k)}(x)\cdot x \A is the inner
     product of \f n^{(k)}(x) \A and \f x \A in \f \RN$. Then the following
     theorem has been obtained in [\JSb] ([\JSb], Theorem 4.6):
\BP

        {\bf Theorem 7.1.} {\it Let \f H \A be as above. Suppose, in 
     addition, that \f \mu \A takes the form of either 
     \f \mu = \mu_0 + \mu_s \A or \f \mu = \mu_0 + \mu_{\ell} $. Let  
     \f \sigma_p(H) \A be the set of the point spectrum of \f H $. Then
     the multiplicity of each \f \lambda \in \sigma_p(H) \A is finite, 
     \f \sigma_p(H) \A does not have any accumulation points except
     at \f 0 \A and \f \infty $.}
\BP

        It is not difficult to extend this result to the general case that
     \f \mu = \mu_0 + \mu_s + \mu_{\ell} $. Using the Example 6.2, we can 
     show a sufficient condition for the nonexistence of the point spectrum 
     of the operator \f H $.
\BP    

        {\bf Theorem 7.2.} {\it Let \f H \A be as above. Suppose that, for 
     almost all \f \omega \in S^{N-1} $, \f \mu_0(r\omega) \A is a
     nondecreasing function of \f r \in [0, \infty) $. Then 
     \f \sigma_p(H) = 0 $.}
\BP

        \demo{ \ \ \ Proof} The condition (ii) of Example 6.2 is now
     satisfied since \f \mu_0(r\omega) \A is nondecreasing. \qed \enddemo
\BP

        Here we are going to give some examples. 
\BP
  
        {\bf Example 7.3.} Let \f \{ U_k \}_{k=0}^{\infty} \A be a sequence
     of open sets of \f \RN \A such that
$$
\left\{ \aligned 
        & \overline{U_k} \subset U_{k+1} \q (k \ge 0), \\
        & \bigcup_{k=0}^{\infty} U_k = \RN, 
\endaligned \right.                                                  \T 7.8
$$
     where the boundary \f \pa U_k \A of \f U_k \A is a continuous surface 
     which is a finite union of smooth surfaces. Suppose that 
$$        
           \tilde n^{(k)}(x)\cdot x \ge 0 \q (k = 0, 1, 2,  \cdots), \T 7.9
$$
     where \f \tilde n^{(k)}(x) \A is the unit outward normal of \f U_k \A
     at \f x \in \pa U_k $. Set
$$
\left\{ \aligned 
        & \Omega_0 = U_0, \\
        & \Omega_k = U_k \bs U_{k-1} \q (k \ge 1), \\
        & S_k^{(+)} = \pa U_k \q \ \ \ \ (k \ge 0), \\
        & S_k^{(-)} = \pa U_{k-1} \q \ (k \ge 1).
\endaligned \right.                                                  \T 7.10
$$
     This is the case that \f K_- = 0 \A and \f K_+ = \infty $. Let 
     \f \mu_0(x) \A be given by
$$
           \mu_0(x) = \nu_k  \qq (x \in \Omega_k),                   \T 7.11
$$
     where \f \{ \nu_k \}_{k=0}^{\infty} \A ia a bounded, positive,
     increasing sequence. Then we see that not only the condition (7.7) is
     satisfied but also \f \mu_0(r\omega) \A is a nondecreasing function of
     \f r \in [0, \infty) \A for almost all \f \omega \in S^{N-1} $.
     Thus Theorem 7.2 can be applied to see that there is no point spectrum
     of \f H $. Therefore the limitig absorption principle holds on the
     whole positive interval \f (0, \infty) \A (see \S5 of [\JSb]).
\BP

        {\bf Example 7.4.} Let \f \{ c_k / k = \pm1, \pm2, \cdots \} \A be
     an increasing sequence of real numbers such that
$$
\left\{ \aligned 
        & c_{-1} < 0 < c_{1}, \\
        & \lim_{k\to\pm\infty} c_k = \pm\infty. \\ 
\endaligned \right.                                                \T 7.12
$$
     Let \f x_N \A be the N-th coordinate of
     \f x = (x_1, x_2, \cdots, x_N) $, and set
$$
      \Omega_k =
\left\{ \aligned 
        & \{ x \in \RN / c_{-1} < x_N < c_{1} \} \qq \ \ \  (k = 0), \\
        & \{ x \in \RN / c_{k-1} < x_N < c_{k} \} 
                                    \qq (k = -1, -2, \cdots), \\ 
        & \{ x \in \RN / c_{k} < x_N < c_{k+1} \} 
                                          \qq (k = 1, 2, \cdots). \\
\endaligned \right.                                                \T 7.13
$$
     We also set

%\newpage

$$
        S_0^{(\pm)} = \{ x \in \RN / x_N = c_{\pm1} \},   \ \ \ \ \ \ \ \ \
                  \ \ \ \ \ \ \ \ \ \ \ \ \ \ \  \ \ \ \ \ \ \ \ \ \ \ \
                                                                    \T 7.14
$$
$$
     S_k^{(+)} =
\left\{ \aligned
        & \{ x \in \RN / x_N = c_{k} \} \qq (k = -1, -2,  \cdots), \\ 
        & \{ x \in \RN / x_N = c_{k+1} \} \q \ \ \ (k = 1, 2, \cdots), \\
\endaligned \right.                                                \T 7.15
$$    
     and
$$
     \ \ \ \ \ \ \ \ \ \ \ \ \ \ \ \ \
        S_k^{(-)} =
\left\{ \aligned
        & \{ x \in \RN / x_N = c_{k-1} \} \q \ \ \ (k = -1, -2,  \cdots), \\ 
        & \{ x \in \RN / x_N = c_{k} \} \qq (k = 1, 2, \cdots), \\
\endaligned \right.                                                  \T 7.16
$$   
     Note that, for \f x \in S_k^{(+)} $, 
$$
     n^{(k)}(x)\cdot x 
\left\{ \aligned
        & \ge 0 \qq (k \ge 0),      \\ 
        & \le 0 \qq (k < 0). \\
\endaligned \right.                                                  \T 7.17
$$       
     Define a simple function \f \mu_0(x) \A by (7.11), where the sequence
     \f \{ \nu_k \}_{k=-\infty}^{\infty} \A is assumed to be bounded and 
     positive such that \f \{ \nu_k \}_{k=-\infty}^{-1} \A is decreasing and 
     \f \{ \nu_k \}_{k=1}^{\infty} \A is increasing. Then, as in Example 7.3,
     Theorem 7.2 can be applied to show that \f \sigma_p(H) = \emptyset $. 
     The planes \f \{ x \in \RN / x_N = c_{k} \} \A can be perturbed as far
     as the condition (7.17) is satisfied.

$$
\ \ \ \
$$

%\newpage

\head   {\bf Appendix}                            \endhead      
\BP
       
        Here we are going to prove a lemma on distributions on a half 
     interval \f (a, \infty) $ which was used when we evaluate the 
     functionals \f M^+(v, r) \A and \f N(v, m, r) $.
\BP

        {\bf Lemma A.} \ {\it Let \f f(r) \A be a real-valued function on 
     \f I = (a, \infty) \A such that \f f \A is locally \f L_1 \A 
     and right continuous on \f I $. Suppose that \f f' \ge 0 $, where
     \f f' \A is the distributional derivative of \f f \A and the 
     inequality should be taken in the sense of distributions. Then
     \f f \A is nondecreasing on \f I $.}
\BP

        \demo{ \ \ \ Proof} Here we are giving a rather elementary proof.
\SP

        (I) Let \f r \in I \A and \f h > 0 $. Then,        
     for \f \phi \in \Con(I) $, we have
$$
\split
       \int_I & [f(r + h) - f(r)]\phi(r) \, dr \\
              & = - \, \int_I f(r)[\phi(r) - \phi(r - h)] \, dr \\
              & = - \, \int_I f(r) \int_{r-h}^r \phi'(s)\, ds \,dr,
\endsplit                                                        \T A.1
$$
     where \f \phi \A is supposed to be extended on the whole line
     \f (-\infty, \infty) \A by setting \f \phi(r) = 0 \A for \f r \le a $.
     Since
$$
           \int_{r-h}^r \phi'(s)\, ds
                    = \int_0^h \phi'(t + r - h) \, dt,             \T A.2
$$
    it follows that
$$
     \int_I [f(r + h) - f(r)]\phi(r) \, dr                            
         = \int_0^h \big[-\, \int_I f(r)\phi'(t + r - h) \, dr \big]\, dt.
                                                                    \T A.3
$$
\SP

        (II) Let \f \phi \in \Con(I) \A and \f \phi \ge 0 $. Then, since
$$
          -\, \int_I f(r)\phi'(t + r - h) \, dr      
                     = <f', \ \phi(\cdot + t - h)> \, \ge \, 0        \T A.4
$$
     for \f h > 0 \A and \f 0 \le t \le h $, where \f <F, \ G> \A
     denotes the value of the distribution \f F \A for the test function
     \f G $, it follows from (A.3) that
$$
             \int_I [f(r + h) - f(r)]\phi(r) \, dr \ge 0         \T A.5
$$
     for any \f \phi \in \Con(I) \A with \f \phi \ge 0 $.   
\SP

        (III) Suppose that there exist \f r_0 \in I $, \f h_0 > 0 \A and
    \f  \eta_0 > 0 \A such that
$$
            f(r_0 + h_0) - f(r_0) = -\, \eta_0.                  \T A.6
$$
     Since \f f \A is right continuous, there exists \f r_1 > r_0 \A
     such that
$$
\left\{ \aligned
          & |f(r_0) - f(r)| < \eta_0/3, \\
          & |f(r_0 + h_0) - f(r + h_0)| < \eta_0/3
\endaligned \right.                                               \T A.7
$$
     for \f r_0 \le r \le r_1 $. Then, for \f r_0 \le r \le r_1 $, we have
$$
\split
      f(r + h_0) & - f(r) \\
                 & = f(r_0 + h_0) - f(r_0)
                              + \{ f(r + h_0) - f(r_0 + h_0) \}
                               + \{ f(r_0) - f(r) \}            \\
                 & \le f(r_0 + h_0) - f(r_0)
                              + | f(r + h_0) - f(r_0 + h_0) |
                                   + | f(r_0) - f(r) |          \\
                 & < -\, \eta_0/3.
\endsplit                                                           \T A.8
$$
     Let \f \phi \in \Con(I) \A such that
$$
\left\{ \aligned
           & \roman{supp}\, \phi \subset [r_0, r_1], \\
           & \phi \ge 0,                              \\
           & \int_{r_0}^{r_1} \phi(r)\, dr = 1.
\endaligned \right.                                               \T A.9
$$
     Then, it follows that
$$
\split
        \int_I [f(r + h) - f(r)]\phi(r) \, dr
                 & = \int_{r_0}^{r_1} [f(r + h) - f(r)]\phi(r) \, dr \\
                 & \le -\,\frac{\eta_0}3\int_{r_0}^{r_1} \phi(r) \, dr \\
                 & = -\,\frac{\eta_0}3 < 0,
\endsplit                                                         \T A.10
$$
     which contradicts (A.5). This completes the proof. \qed \enddemo         

$$
\ \ \
$$

%\newpage

\head {\bf Refenreces}                    \endhead
\BP

\ref \by [\Eia] D. Eidus \paper The principle of limiting absorption \yr 1965 
\vol 47 \jour American Math. Soc. Translations
\pages 157-191 (Mat. Sb. {\bf 57} (1962))
\endref
\SP

\ref \by [\Eib] D. Eidus \pages 29-38 \paper The limiting absorption and
amplitude problems for the diffraction problem with two unbounded media 
\yr 1986 \vol 107 \jour Comm. Math. Phys.
\endref
\SP

\ref \by [\Ik] T. Ikebe \paper Eigenfunction expansions associated with
the \Sch\ operators and their applications to scattering theory \jour Arch.
Rational Mech. Anal. \vol 5 \yr 1960 \pages 1-34
\endref
\SP

\ref \by [\IS] T. Ikebe and Y. Sait\B o \paper Limiting absorption method and
absolute continuity for the Schr\"odinger operator \jour J. Math. Kyoto Univ.
\vol 12 \yr 1972 \pages 513-612.
\endref
\SP

\ref \by [\IU] T. Ikebe and J. Uchiyama \paper On the asymptotic behavior
of eigenfunctions of second-order elliptic differential operators \jour J.
Math. Kyoto Univ. \vol 11 \yr 1971 \pages 425-448  
\endref
\SP

\ref \by [\Ja] W. J\"ager \paper Zur Theorie der Schwingugsgleichung mit 
variablen Koeffizienten in Aussengebieten \jour Math. Z. \vol 102 \yr 1969 
\pages 62-88 
\endref
\SP

\ref \by [\JSa] W. J\"ager and Y. Sait\B o \paper On the Spectrum of the 
Reduced Wave Operator with Cylindrical Discontinuity \jour Forum
Mathematicum \vol 9 \yr 1997 \pages 29-60
\endref
\SP

\ref \by [\JSb] W. J\"ager and Y. Sait\B o \paper The reduced wave equation 
in layered materials \jour to appear in Osaka J. Math.  
\endref
\SP
      
\ref \by [\Ka] T. Kato \paper Growth properties of solutions of the reduced
wave equation with a variable coefficients \jour Comm. Pure Appl. Math.
\vol 12 \yr 1959 \pages 403-425  
\endref
\SP

\ref \by [\Mi] C. M\"uller \book Grundprobleme der mathematischen Theorie
elektromagnetischer \linebreak Schwingungen \publ Springer, Berlin \yr 1957  
\endref
\SP
   
\ref \by [\Rl] F. Rellich \paper \"Uber das asymptotische Verhalten des
L\"osungen von \f \Delta u + \lambda u = 0 $ \jour Jber. Deutsche. Math.
Verein. \vol 53 \yr 1943 \pages 57-65  
\endref
\SP
   
\ref \by [\RZ] G. Roach and B. Zhang \paper On Sommerfeld radiation conditions for 
the diffraction problem with two unbounded media \jour Proc. Royal Soc. 
Edinburgh \vol 121A \yr 1992 \pages 149-161 
\endref
\SP

\ref \by [\YS] Y. Sait\B o \book Spectral Representations for Schr\"odinger 
Operators with Long-Range Potentials, Lecture Notes in Mathematics 
\publ Springer, Berlin \vol 727 \yr 1979  
\endref
\SP
 
\ref \by [\We] J. Weidmann \paper On the Continuous spectrum of \Sch\ operators
\jour Comm. Pure and Appl. Math. \vol 19 \yr 1966 \pages 107-110
\endref
\SP

\ref \by [\Wi] E. Wienholtz \paper Halbbeschr\"anke partielle 
Differentialoperatoren zweiter Ordnung vom elliptischen Typus \jour Math.
Ann. \vol 135 \yr 1958 \pages 50-80 
\endref
\SP

%\ref \by  \paper \jour \publ  \vol  \yr \pages  
%\endref
%\SP
   
\enddocument